\documentclass[12pt]{amsart}%
\usepackage{amsfonts}
\usepackage{latexsym}
\usepackage{amssymb}
\usepackage{amsmath}
\usepackage{amstext}
\usepackage{amsthm}
\usepackage{amsxtra}
\usepackage[utf8]{inputenc}
\usepackage[brazilian,english]{babel}
\providecommand{\U}[1]{\protect\rule{.1in}{.1in}}
\newtheorem{theorem}{Theorem}
\newtheorem{lemma}{Lemma}

\newtheorem{example}{Example}
\newtheorem{proposition}{Proposition}
\newtheorem{remark}{Remark}
\theoremstyle{definition}

\usepackage{enumerate}
\newtheorem*{theoremA}{Theorem A}
\newtheorem*{theoremB}{Theorem B}

\def\D{{\mathbb D}}

\def\({\left(}       \def\){\right)}

\begin{document}
	
	\title[Spaces of harmonic mappings and the Schwarzian derivative]{Properties of Besov and $Q_p$ spaces in terms of the Schwarzian derivative of harmonic mappings}
	
	\author[H. Arbel\'aez]{Hugo Arbel\'{a}ez}
	\address{Escuela de Matem\'aticas, Universidad Nacional de Colombia, Medell\'{\i}n, Colombia}
	\email{hjarbela@unal.edu.co}
	
	\author[R. Hern\'andez]{Rodrigo Hern\'andez}
	\address{Facultad de Ingenier\'{\i}a y Ciencias, Universidad Adolfo Ib\'a\~nez. Av. Pa\-dre Hurtado 750. Vi\~na del Mar, Chile.} \email{rodrigo.hernandez@uai.cl}
	
	\author[W. Sierra]{Willy Sierra}
	\address{Departamento de Matem\'aticas, Universidad del Cauca, Popay\'{a}n, Colombia}
	\email{wsierra@unicauca.edu.co}
	

	\thanks{The first author was partially supported by the Universidad Nacional de Colombia, Hermes code 61126. The second author is supported by grant Fondecyt $\#1190756$, Chile. The third author thanks the Universidad del Cauca for providing time for this work through research project VRI ID 6235.}

	\maketitle

	\begin{abstract}
		In this paper we give a characterization of $\log J_f$ belongs to $\widetilde{\mathcal{B}}_p$ or $\widetilde{\mathcal{Q}}_p$ spaces for any locally univalent sense-preserving harmonic mappings $f$ defined in the unit disk, using the Schwarzian derivative of $f$ and Carleson meseaure. In addition, we introduce the classes $\mathcal{BT}_p$ and $\mathcal{QT}_p$, based on the Jacobian operator, and begin a study of these.
	\end{abstract}
	
	
	
	\section{Introduction and preliminaries}\label{Section preliminaries}

	The Besov spaces $\mathcal{B}_p$ have been extensively studied since their introduction by Zhu in \cite{Zhu}. These spaces can be seen as a particular case of weighted Bergman spaces, defined as the analytic functions on the unit disk whose derivative is integrable with respect to the weight measure $(1-|z|^2)^\alpha$, in this case, $\alpha = p-2$. Readers can find more information on Bergman spaces in the excellent book by P. Duren and A. Schuster \cite{DS}. On the other hand, the $Q_p$ spaces introduced by J. Xiao in \cite{Xiao} are also defined for those analytic functions $f$ such that $\sup{\int_\D|f'(z)|^2(1-|\sigma_a(z)|^2)^p dA(z)}$ is finite. In this brief manuscript, we will focus on studying the relationship of these spaces with the Schwarzian derivative and Carleson measures, in the context of complex harmonic functions that preserve orientation. More specifically, we will extend the following results that account for these relationships in the realm of analytic functions:
	
	\begin{theoremA} Let $f:\D\to\Omega$ be a conformal map such that $f(\partial \D)$ is a closed Jordan curve. Then
		\begin{itemize}
			\item[i)] If $1<p<\infty,$ $\log f'$ is in $\mathcal{B}_p$ if and only if $\displaystyle \int_{\D}|Sf(z)|^p(1-|z|^2)^{2p-2}dA(z)<\infty$.\\
			
			\item[ii)] If $0<p\leq 1,$ $\log f'$ is in $Q_{p,0}$ if and only if $|Sf(z)|^2(1-|z|^2)^{2+p}dA(z)$ is a vanishing $p$-Carleson measure.
		\end{itemize}
	\end{theoremA}
	
	\begin{theoremB} Let $0<p<\infty$ and $f:\D\to\Omega$ be a conformal map, then $\log f'\in Q_p$ if and only if $|Sf(z)|^2(1-|z|^2)^ {2+p}dA(z)$ is an $p$-Carleson measure.
	\end{theoremB}
	
	These results appear in \cite{Fernando Jouni} and \cite{Pau Pelaez}, respectively.
	
	In Subsection\,\ref{subsection harmonic spaces}, we focus on the definition of the spaces of Besov $\widetilde{\mathcal{B}}_p$ and $\widetilde{\mathcal{Q}}_p$ for smooth mappings $F$ such that the derivatives with respect to $z$ and $\overline{z}$ exist in $\D$. Based on these definitions, in Section\,\ref{section main results} our efforts are dedicated to proving Theorems A and B for these types of spaces. Finally, in Section\,\ref{section results in harmonic spaces} functions called Besov type and $\widetilde{\mathcal{Q}}_p$ type are introduced in a similar manner to what is done in \cite{Bloch type map}, and we note, in Remark\,\ref{rem 4}, the difficulties that these types of functions present. 
	
	As we have already mentioned, in these brief notes we address the aforementioned results but in the context of locally univalent complex harmonic functions defined on the unit disk. Using as a principal tools the Schwarzian derivative defined in \cite{HM} for this kinds of mappings, and we defined the corresponding Besov and $\widetilde{\mathcal{Q}}_p$ spaces for this scenario. 
	
	\subsection{Planar harmonic mappings}
	
	In this section we introduce some notation and we present several classical results concerning harmonic mappings in the plane. We refer to \cite{Du04} and to the references therein for more details about this topic.
	
	Let $f$ be a planar harmonic mapping defined on the unit disk $\mathbb{D}=\{z: |z|<1\}.$ In this case, $f$ has the canonical representation $f=h+\overline{g}$, where $h$ and $g$ are analytic functions in $\mathbb{D}$; the representation is unique under the condition that $g(z_0)=0$ for some $z_0$ fixed in $\mathbb{D}$. Here we will assume the condition $h(0)=g(0)=0.$ It is well known that $f$ is locally univalent if and only if its Jacobian $J_f=|h'|^2-|g'|^2$ does not vanish. Thus, if $f$ is locally univalent, it is either sense-preserving or sense-reversing depending on the conditions $J_f>0$ or $J_f<0$ throughout its domain, respectively.  Along this paper we will consider sense-preserving harmonic mappings on $\mathbb{D},$ case in which the analytic part $h$ is locally univalent in $\mathbb{D}$ and the second complex dilatation of $f,$ $w = g'/h',$ is an analytic function in $\mathbb{D}$ satisfying the condition $|w(z)|<1,$ for all $z\in\mathbb{D}.$ 
	
	We define the function $F=\log J_f.$ By using the operators $\partial_z$ and $\partial_{\overline{z}},$ defined by
	\[\partial_z:=\frac{1}{2}\left(\partial_x-i\partial_y \right)\qquad\text{and}\qquad \partial_{\overline{z}}:= \frac{1}{2}\left(\partial_x+i\partial_y \right),\]
	we have $F_{\overline{z}}=\overline{F_z},$
	\begin{equation}\label{F sub z and F sub zz}
		F_z=\frac{h''}{h'}-\frac{w'\bar{w}}{1-|w|^2},\qquad\text{and}\qquad F_{zz}=\left( \frac{h''}{h'}\right)'-\frac{w''\bar{w}}{1-|w|^2}-\left(\frac{w'\bar{w}}{1-|w|^2} \right)^2.
	\end{equation}
	Note that $F_z$ coincides with the Pre-Schwarzian derivative of the harmonic function $f,$ denoted by $P_f,$ which was introduced in \cite{HM}. In that same paper the authors define the Schwarzian derivative of $f$ by $S_f=\partial_z P_{f}-\frac{1}{2}(P_f)^2.$ It follows from \eqref{F sub z and F sub zz} that
	\begin{equation}\label{Schwarzian of f}
		\begin{split}
			S_f&=F_{zz}-\frac{1}{2}(F_z)^2\\
			&=S_h+\frac{\bar{w}}{1-|w|^2}\left(w'\frac{h''}{h'}-w'' \right)
			-\frac{3}{2}\left(\frac{w'\bar{w}}{1-|w|^2} \right)^2,
		\end{split}
	\end{equation}
	where $S_h$ denotes the Schwarzian derivative of the analytic function $h,$ which is given by
	\[S_h=\left(\frac{h''}{h'} \right)'-\frac{1}{2}\left(\frac{h''}{h'} \right)^2. \]
	As in the analytic case, the Pre-Schwarzian and Schwarzian norms of $f$ are given by
	\[\|P_f\|:=\sup_{z\in\mathbb{D}}\{(1-|z|^2)|P_f(z)|\}\qquad\text{and}\qquad \|S_f\|:=\sup_{z\in\mathbb{D}}\{(1-|z|^2)^2|S_f(z)|\},  \]
	respectively. Both in the analytic case and in the harmonic case these norms are closely related with the notion of uniformly locally univalent functions on the unit disk, that is, those harmonic functions $f:\mathbb{D}\to\mathbb{C}$ for which there is a constant $0<\rho:=\rho(f)\leq \infty,$ such that $f$ is univalent in each hyperbolic disk with center $a\in\mathbb{D}$ and hyperbolic radius $\rho,$ which will be denoted by $D(a,\rho).$ Note that $f$ is univalent if $\rho=\infty.$ We remember that the hyperbolic distance in $\mathbb{D}$ is defined by
	\[ d_h(z,\xi)=\frac{1}{2}\log\frac{1+p(z,\xi)}{1-p(z,\xi)},\qquad\text{where}\quad p(z,\xi)=\left|\frac{z-\xi}{1-\overline{\xi}z}\right|. \]
	When $f$ is analytic, in \cite{ULU Yamashita} is proven that $f$ is uniformly locally univalent if and only if $\|P_f\|$ is finite, which also is equivalent to $\|S_f\|$ finite. A generalization to the harmonic case can be found in \cite{HM}, where it is also proven that for a locally univalent harmonic mapping $f=h+\overline{g},$ $\|S_f\|$ is finite if and only if $\|S_h\|$ is finite; we will use these facts later without further comment.
	
	\subsection{Some spaces of harmonic mappings}\label{subsection harmonic spaces} Now we define the functions spaces that we will use in this paper. We start recalling that if $A(\mathbb{D})$ denotes the space of analytic functions in $\mathbb{D},$ the analytic Bloch space is defined by
	\[\mathcal{B}=\{h\in A(\mathbb{D}): \|h\|_{\mathcal{B}}:=\sup_{z\in\mathbb{D}}(1-|z|^2)|h'(z)|<\infty \}\]
	and the little analytic Bloch space is the subspace of $\mathcal{B}$ given by
	\[\mathcal{B}_0:=\{h\in A(\mathbb{D}): \lim_{|z|\to 1^{-}}(1-|z|^2)|h'(z)|=0 \}.\]
	For $1<p<\infty,$ we say that $h\in A(\mathbb{D})$ belongs to the analytic Besov space $\mathcal{B}_p$ if
	\begin{equation}\label{def norm Besov space}
		\|h\|_{\mathcal{B}_p}^p:=\int_{\mathbb{D}}|h'(z)|^p(1-|z|^2)^{p-2}dA(z)<\infty,
	\end{equation}
	where $dA$ is the element of the Lebesgue area measure on $\mathbb{D}.$ The space $Q_p,$ $0<p<\infty,$ consists of all functions $h\in A(\mathbb{D})$ satisfying the condition
	\begin{equation*}
		\sup_{a\in\mathbb{D}}\int_{\mathbb{D}}|h'(z)|^2(1-|\varphi_a(z)|^2)^{p}dA(z)<\infty,
	\end{equation*}
	where $\varphi_a(z):=\frac{a-z}{1-\overline{a}z}.$ Note that $\varphi_a$ is an automorphism of $\mathbb{D}$ satisfying $\varphi_a(a)=0,$ $\varphi_a(0)=a,$ $(1-|z|^2)|\varphi_a'(z)|=1-|\varphi_a(z)|^2,$ and $\varphi_a^{-1}=\varphi_a.$ We will say that a function $h\in A(\mathbb{D})$ belongs to the space $Q_{p,0}$ if
	\begin{equation*}
		\lim_{|a|\to1^-}\int_{\mathbb{D}}|h'(z)|^2(1-|\varphi_a(z)|^2)^{p}dA(z)=0.
	\end{equation*}
	Throughout the paper we use $a\lesssim b$ $(a\gtrsim b)$ for denote $a\leq C b$ $ (a\geq Cb),$ being $C$ a constant independent of $a$ and $b.$ We also write $a\simeq b,$ if $a\lesssim b$ and $a\gtrsim b.$ With this notation we have the following results, which will be a key ingredient in
	the proofs of the proposed results. The first of them, see for example \cite{Pau Pelaez}, establishes that if $0<p<\infty$ and $-1<\alpha<\infty,$
	\begin{equation}\label{eq. norm bergman spaces}
		\int_{\mathbb{D}}|h(z)|^p(1-|z|^2)^{\alpha}dA(z)\simeq \int_{\mathbb{D}}|h'(z)|^p(1-|z|^2)^{p+\alpha}dA(z)+|h(0)|^p, 
	\end{equation}
	for all $h\in A(\mathbb{D}).$ The second, which is proven in \cite[Theorem 1]{Qs and n derivative}, says that for $n\geq 1$ and $0<p<\infty,$ $h\in Q_p$ if and only if
	\begin{equation}\label{eq. f in Qp iff ...}
		\sup_{a\in\mathbb{D}}\int_{\mathbb{D}}|h^{n}(z)|^2\left(1-|\varphi_a(z)|^2 \right)^p(1-|z|^2)^{2n-2}dA(z)<\infty.
	\end{equation} 
	The classical theory of analytic functions spaces has a natural extension to the setting  of smooth functions from $\mathbb{D}$ into $\mathbb{C},$ which has been being investigated in the last years. For some works on this topic, we refer to the reader to \cite{Thesis Aljuaid,Colona 1989 def Bloch,Bloch type map,Yoneda Bloch Besov Harmonic}. In particular, in \cite{Colona 1989 def Bloch} is studied the harmonic extension of $\mathcal{B}$ and it is shown that a harmonic function $F:\mathbb{D}\to\mathbb{C}$ belongs to the hamonic Bloch space if and only if
	\begin{equation}\label{eq. def H-Bloch}
		\sup_{z\in\mathbb{D}}\{(1-|z|^2)(|F_z(z)|+|F_{\bar{z}}(z)|)\}<\infty.
	\end{equation}
	This definition clearly can be extended to include the whole family of smooth functions from $\mathbb{D}$ into $\mathbb{C}.$ Thus, a smooth function $F:\mathbb{D}\to\mathbb{C}$ is said to be a Bloch function if it satisfies \eqref{eq. def H-Bloch}. We will use $\widetilde{\mathcal{B}}$ to denote the space of such Bloch functions.
	
	In a similar way to the previous case, we use the definition proposed in \cite{Thesis Aljuaid} of harmonic Besov spaces to  include smooth functions from $\mathbb{D}$ into $\mathbb{C}.$ For $1<p<\infty,$ if a smooth function $F:\mathbb{D}\to\mathbb{C}$ satisfies
	\begin{equation}\label{eq. def H-Besov}
		\int_{\mathbb{D}}(|F_z(z)|+|F_{\bar{z}}(z)|)^{p}(1-|z|^2)^{p-2}dA(z)<\infty,
	\end{equation}
	we will say that it belongs to the Besov space $\widetilde{\mathcal{B}}_p.$
	\begin{remark}\label{rem log h' and log Jf}
		If $f=h+\overline{g}$ is a harmonic function and $F=\log J_f,$ then \eqref{F sub z and F sub zz} gives us
		\begin{equation*}
			\left| \frac{h''}{h'} \right|^p\lesssim |F_z|^p + \left( \frac{|w'\overline{w}|}{1-|w|^2} \right)^p\quad\text{and}\quad |F_z|^p\lesssim \left| \frac{h''}{h'} \right|^p + \left( \frac{|w'\overline{w}|}{1-|w|^2} \right)^p,
		\end{equation*}
		$1<p<\infty.$ Hence, if $\left\| w\right\|:=\sup \{|w(z)|: z\in\mathbb{D}\} <1$ and $w\in\mathcal{B}_p,$ then $F\in\widetilde{\mathcal{B}}_p$ if and only if $\log h'\in\mathcal{B}_p.$
	\end{remark}
	By following the same line as in \cite{Qp H-spaces}, a smooth function $F:\mathbb{D}\to \mathbb{C}$ belongs to the space $\widetilde{\mathcal{Q}}_p,$ $0<p<\infty,$ if it satisfies
	\begin{equation}\label{eq. def of Qs}
		\sup_{a\in\mathbb{D}}\int_{\mathbb{D}}(|F_z(z)|+|F_{\bar{z}}|)^2(1-|\varphi_a(z)|^2)^{p}dA(z)<\infty
	\end{equation}
	and $F\in\widetilde{\mathcal{Q}}_{p,0},$ if
	\begin{equation*}
		\lim_{|a|\to 1^-}\int_{\mathbb{D}}(|F_z(z)|+|F_{\bar{z}}|)^2(1-|\varphi_a(z)|^2)^{p}dA(z)=0.
	\end{equation*}
	\subsection{Carleson measures} With a similar notation to that used in \cite{Pau Pelaez,Fernando Jouni}, given an arc $I\subseteq\partial\mathbb{D},$ we write $|I|$ for the normalized arclength  of $I$ and we define the set
	\[ S(I):=\{z=re^{it}\in\mathbb{D}: 1-r<|I|\;\text{ and }\; e^{it}\in I \} .\] 
	Let $\mu$ be a positive Borel measure on $\mathbb{D}$ and $p>0.$ $\mu$ is said to be a $p-$Carleson measure on $\mathbb{D},$ if
	\[\|\mu\|_p:=\sup_{I\subseteq\partial\mathbb{D}}\frac{\mu(S(I))}{|I|^p}<\infty,\]
	which is equivalent to the condition
	\begin{equation}\label{eq. p carleson measure}
		\sup_{a\in\mathbb{D}}\int_{\mathbb{D}}\left|\varphi_a'(z) \right|^pd\mu(z)<\infty.
	\end{equation}
	A $p-$Carleson measure $\mu$ on $\mathbb{D}$ is called a $p$-vanishing Carleson measure on $\mathbb{D}$ if
	\begin{equation*}
		\lim_{|I|\to 0}\frac{\mu(S(I))}{|I|^p}=0.
	\end{equation*}
	This is equivalent to the condition
	\begin{equation*}
		\lim_{|a|\to1^-}\int_{\mathbb{D}}\left|\varphi_a'(z) \right|^pd\mu(z)=0, 
	\end{equation*}
	see \cite[Lemma 2.1]{BMOA Carleson}. We are now ready to address the main theme of the paper.
	\section{Relationship between Besov spaces $\widetilde{\mathcal{B}}_p$ and $\widetilde{\mathcal{Q}}_p$ spaces and the Schwarzian derivative}\label{section main results}
	
	In this section we will study in the setting of harmonic mappings some of the results presented in \cite{Fernando Jouni} for the case of analytic functions. We start with the following theorem, which is a harmonic version of one of the implications of part i) of Theorem\,A (Theorem\,1 in \cite{Fernando Jouni}).
	
	\begin{theorem}\label{thm. Jf implies I finite}
		Let $1<p<\infty$ and $f:\mathbb{D}\to\mathbb{C}$ be a uniformly locally univalent harmonic map with dilatation $w\in \mathcal{B}_p$ satisfying $\|w\|<1.$ If $F=\log J_f\in \widetilde{\mathcal{B}}_p,$ then
		\begin{equation}\label{eq. I de f}
			I(f):=\int_{\mathbb{D}}|S_f(z)|^p(1-|z|^2)^{2p-2}dA(z)
		\end{equation}
		is finite.
	\end{theorem}
	\begin{proof}
		We suppose that $F=\log J_f\in \widetilde{\mathcal{B}}_p$ and we write $f$ in the standard form $h+\overline{g}.$  By equation \eqref{Schwarzian of f} one has that
		\begin{equation}\label{est I(f)}
			I(f)\lesssim \int_{\mathbb{D}}|F_{zz}(z)|^p(1-|z|^2)^{2p-2}dA(z)+\int_{\mathbb{D}}|F_z(z)|^{2p}(1-|z|^2)^{2p-2}dA(z).
		\end{equation}
		On the other hand, the condition $\|w\|<1$ and \eqref{eq. norm bergman spaces} imply
		\begin{equation}\label{Aux 1}
			\begin{split}
				\int_{\mathbb{D}}(1-|z|^2)^{2p-2}\left| \frac{w''(z)\overline{w(z)}}{1-|w(z)|^2} \right|^pdA(z)&\lesssim \int_{\mathbb{D}}\left| w''(z) \right|^p(1-|z|^2)^{p+(p-2)}dA(z)\\
				&\lesssim \int_{\mathbb{D}}|w'(z)|^p(1-|z|^2)^{p-2}dA(z)-|w(0)|^p,
			\end{split}
		\end{equation}
		which is finite since $w\in \mathcal{B}_p.$
		
		Next, from $\|w\|<1$ we have the estimate
		\begin{equation*}
			(1-|z|^2)^{2p-2}\left| \frac{w'(z)\overline{w(z)}}{1-|w(z)|^2} \right|^{2p}\lesssim (1-|z|^2)^{p-2}|w'(z)|^p\left| \frac{(1-|z|^2)w'(z)}{1-|w(z)|^2} \right|^{p}.
		\end{equation*}
		It follows from the Schwarz-Pick inequality and $w\in \mathcal{B}_p,$ that
		\begin{equation}\label{Aux 2}
			\int_{\mathbb{D}}(1-|z|^2)^{2p-2}\left| \frac{w'(z)\overline{w(z)}}{1-|w(z)|^2} \right|^{2p} dA(z)<\infty.
		\end{equation}
		On the other hand, we note that in virtue of \eqref{eq. norm bergman spaces},
		\begin{equation*}
			\begin{split}
				\int_{\mathbb{D}}(1-|z|^2)^{2p-2}\left| \left(\frac{h''(z)}{h'(z)} \right)' \right|^pdA(z)&= \int_{\mathbb{D}}(1-|z|^2)^{p+(p-2)}\left| \left(\frac{h''(z)}{h'(z)} \right)' \right|^pdA(z)\\
				&\lesssim \int_{\mathbb{D}}\left|\frac{h''(z)}{h'(z)} \right| ^p(1-|z|^2)^{p-2}dA-\left| \frac{h''(0)}{h'(0)}\right| ^p<\infty,
			\end{split}
		\end{equation*}
		since by hypothesis and Remark\,\ref{rem log h' and log Jf}, $\log h'\in\mathcal{B}_p.$ We conclude from this, \eqref{F sub z and F sub zz}, \eqref{Aux 1}, and \eqref{Aux 2}, that the first integral in \eqref{est I(f)} is finite.
		
		With respect to the second integral of \eqref{est I(f)}, we observe that
		\begin{equation}\label{eq. Aux 3}
			|F_z(z)|^{2p}(1-|z|^2)^{2p-2}=\left( |F_z(z)|(1-|z|^2)\right)^p |F_z(z)|^{p}(1-|z|^2)^{p-2}.
		\end{equation}
		Now, the fact that $f$ is uniformly locally univalent  guarantees us that $\|S_f\|<\infty$ and therefore $\|S_h\|<\infty$ \cite[Theorem 6]{HM}.  Hence, $h$ is uniformly locally univalent, or equivalently, $\|P_h\|<\infty.$ It follows from \eqref{F sub z and F sub zz} and the Schwarz-Pick inequality that the first factor of the right side of \eqref{eq. Aux 3} is bounded. In consequence,
		\[\int_{\mathbb{D}}|F_z(z)|^{2p}(1-|z|^2)^{2p-2}dA(z)\lesssim \int_{\mathbb{D}} |F_z(z)|^{p}(1-|z|^2)^{p-2} dA(z)<\infty, \]
		since $F\in\widetilde{\mathcal{B}}_p.$ We conclude that $I(f)<\infty,$ which ends the proof.
	\end{proof}
	\begin{remark}\label{rem. I(h) I(f) finite}
		Following similar arguments to that given in the above proof, one can show that under the same conditions on $w,$ $I(f)$ is finite if and only if $I(h)$ is finite. We will use this fact to give other alternative approximation of part i) of Theorem\,A. To this end we will consider a linear combination of the form $\varphi_{\lambda}:=h+\lambda g,$ $\lambda\in\partial\mathbb{D},$ where $f=h+\overline{g}.$ We remark that the relation between the analytic function $\varphi_{\lambda}$ and the harmonic function $f$ has been used by several authors to obtain properties of $f$ from properties of $\varphi_{\lambda}.$ An important example of this fact is the called shear construction, which was introduced in \cite{Clunie and Sheil} to construct sense preserving univalent harmonic mappings in the unit disk. See also \cite{HeMa13} for other results relating properties of $\varphi_{\lambda}$ with those of the corresponding function $f.$   
	\end{remark}
	\begin{theorem}\label{Thm 2}
		Let $1<p<\infty$ and $f=h+\overline{g}:\mathbb{D}\to\mathbb{C}$ be a univalent harmonic map such that its dilatation $w\in \mathcal{B}_p$ satisfies $\|w\|<1$ and assume that $\varphi:=\varphi_{\lambda}=h+\lambda g$ is a conformal map with $\varphi(\partial\mathbb{D})$ a closed Jordan curve for some $\lambda\in\partial\mathbb{D}.$ Then $F=\log J_f\in \widetilde{\mathcal{B}}_p$ if and only if $I(f)$ is finite.
	\end{theorem}
	\begin{proof}
		That $F=\log J_f\in \widetilde{\mathcal{B}}_p$ implies $I(f)<\infty,$ is the statement of the above theorem. To the converse, from the definition of $\varphi,$ we obtain
		\[\varphi'=h'(1+\lambda w)\qquad\text{and}\qquad \frac{\varphi''}{\varphi'}=\frac{h''}{h'}+\frac{\lambda w'}{1+\lambda w},\]
		from which we get, after a straightforward calculation, that
		\begin{equation}\label{eq. Sh and S(phi)}
			S_h=S_{\varphi}+\frac{\varphi''}{\varphi'}\frac{\lambda w'}{1+\lambda w}+\frac{1}{2}\left(\frac{\lambda w'}{1+\lambda w} \right)^2-\frac{\lambda w''}{1+\lambda w}. 
		\end{equation}
		We conclude, following the same arguments as in the proof of Theorem\,\ref{thm. Jf implies I finite}, that $I(h)<\infty$ if and only if $I(\varphi)<\infty.$ Since Remark\,\ref{rem. I(h) I(f) finite} and the hypothesis imply $I(h)<\infty,$ it follows that $I(\varphi)<\infty.$ Now we use Theorem\,1 in \cite{Fernando Jouni} to obtain that $\log\varphi'\in \mathcal{B}_p,$ which implies that
		\[\log h'=\log\varphi'-\log(1+\lambda w)\in \mathcal{B}_p.\]
		The theorem follows by using Remark\,\ref{rem log h' and log Jf}.
	\end{proof}
	\begin{example}\label{exa shear construction}
		To illustrate Theorem\,\ref{Thm 2}, we consider the shear $f=h+\bar{g}$ of the identity function $\varphi(z)=z$ with dilatation $w(z)=\rho z,$ $0<\rho<1.$ Then $h-g=\varphi$ and $g'=w h',$ whence
		\[h'(z)=\frac{1}{1-\rho z}\qquad\text{and}\qquad g'(z)=\frac{\rho z}{1-\rho z}.\]
		Moreover, it is clear that $\|w\|<1$ and $w\in \mathcal{B}_p,$ for all $p>1.$ Note that $\log h'\in \mathcal{B}_p$ and therefore Remark\,\ref{rem log h' and log Jf} implies that $\log J_f\in \widetilde{\mathcal{B}}_p,$ for all $p>1.$ We conclude from Theorem\,\ref{Thm 2} that $I(f)<\infty.$ 
	\end{example}
	\begin{remark}\label{rem. norm of w}
		In relation with the assumption $\|w\|<1,$ we highlight that this condition is not so restrictive; it in general can not be omitted from the statement of the above results. For example, we will show that most of harmonic function $f=h+\bar{g},$ with dilatation $w(z)=z,$ satisfy $\log J_f\notin \widetilde{\mathcal{B}}_p,$ for all $p.$ Indeed, if $w(z)=z$ one has that
		\begin{equation*}
			\log J_f\in\widetilde{\mathcal{B}}_p\quad\text{if and only if}\quad\int_{\mathbb{D}} \left|\partial_z\log J_f(z) \right|^p(1-|z|^2)^{p-2}dA(z)<\infty,
		\end{equation*}
		whence $\log J_f\in\widetilde{\mathcal{B}}_p$ if and only if
		\begin{equation}\label{eq. aux remark 3}
			\int_{\mathbb{D}} \left|(1-|z|^2)\frac{h''(z)}{h'(z)}-\bar{z} \right|^p(1-|z|^2)^{-2}dA(z)<\infty.   
		\end{equation}
		However, assuming for example that $h''/h'$ has finite angular limit at some point $\xi\in\partial\mathbb{D}$ (this is essentially the case for most of meromorphic functions in $\mathbb{D},$ except possibly those with a behaviour very bad at almost all points of $\partial\mathbb{D},$ see Plessner's Theorem in \cite{Book_Pommerenke}), it can be proven that the integral in \eqref{eq. aux remark 3} diverges. In effect, in this case we can choose a Stolz angle at $\xi$ of the form
		\[ \Delta= \{ z\in\mathbb{D} : |\arg(1-\overline{\xi}z)|<\alpha\text{ and } |z-\xi|<\rho<1/2  \},\]
		with $0<\alpha<\frac{\pi}{2}$, $0<\rho<2\cos\alpha,$ and $0<\delta_0<2\rho /3$ such that
		\begin{equation*}
			\left|(1-|z|^2)\frac{h''(z)}{h'(z)}-\bar{z} \right|\geq |z|-\frac{1}{2}\geq (1-\rho)-\frac{1}{2}=\frac{1}{2}-\rho,
		\end{equation*}
		for all $z\in \Delta\cap B(\xi,\delta_0),$ where $B(a,r)$ denote the euclidean ball of radius $r>0$ and center $a.$  Thus, for all $0<\delta<\delta_0,$
		\begin{align*}
			\int_{\Omega_{\delta}} \left|(1-|z|^2)\frac{h''(z)}{h'(z)}-\bar{z} \right|^p(1-|z|^2)^{-2}dA(z)&\geq \left(\frac{1}{2}-\rho\right)^p\int_{\Omega_{\delta}}\frac{dA(z)}{(1-|z|^2)^2}\\
			&\geq \left(\frac{1}{2}-\rho\right)^p\frac{1}{\left(1-(1-\delta)^2\right)^2}\int_{\Omega_{\delta}}dA(z),
		\end{align*}
		where $\Omega_{\delta}:=\Delta\cap B(\xi,\delta).$ In consequence,
		\[\int_{\Omega_{\delta}} \left|(1-|z|^2)\frac{h''(z)}{h'(z)}-\bar{z} \right|^p(1-|z|^2)^{-2}dA(z)\geq C\frac{\delta^2}{4\delta^2-4\delta^3+\delta^4},\]
		for some constant $C$ independent of $\delta.$ By taking limit when $\delta\to 0$ it follows that the integral in \eqref{eq. aux remark 3} must be divergent and therefore $\log J_f\notin \widetilde{\mathcal{B}}_p,$ for all $p.$
		
		The same argument shows that most of harmonic function $f=h+\bar{g},$ with dilatation an automorphism of the disk, satisfy $\log J_f\notin \widetilde{\mathcal{B}}_p,$ for all $p.$   
	\end{remark}
	Next, we will prove results similar to the previous ones but this time in spaces $\widetilde{\mathcal{Q}}_p.$ To this end we need a version in the setting of uniformly locally univalent analytic functions of one of the implications of Theorem\,B; more precisely we require the following result:
	\begin{proposition}\label{Prop Thm of Pau-Pelaez for locally univalent}
		Suppose that $0<p<\infty$ and $h$ is a uniformly locally univalent analytic function in the unit disk. If $|S_h(z)|^2(1-|z|^2)^{2+p}$ is a $p-$Carleson measure, then $\log h'\in Q_p.$ 
	\end{proposition}
	The proof of the proposition is a slight modification of the proof of  (b) implies (a) in \cite[Theorem 1]{Pau Pelaez}  and we will omit it. We remark that except for the value of the constants, the previous results required for the proof of (b) implies (a) in \cite[Theorem 1]{Pau Pelaez} are also true for  uniformly locally univalent analytic function. For example, for functions $h$ in this class there is $M>0$ such that $(1-|z|^2)^2|S_h(z)|\leq M,$ for all $z\in\mathbb{D}$ (compare with Lemma\,A in \cite{Pau Pelaez}). For the sake of completeness we prove in this context a version of Lemma\,4 in \cite{Pau Pelaez}.   
	\begin{lemma}
		Let $h$ be a uniformly locally univalent analytic function in $\mathbb{D}$ and suppose that there is $z_0\in\mathbb{D}$ such that
		\[(1-|z_0|^2)^2|S_h(z_0)|>\delta.\]
		Then there is a positive constant $c=c(\delta,h)<1$ such that
		\[(1-|z|^2)^2|S_h(z)|>\frac{\delta}{32},\]
		for all $z\in B(z_0, c(1-|z_0|^2)).$
	\end{lemma}
	\begin{proof}
		Let $\rho>0$ such that $h$ is univalent in each hyperbolic disk $D(a,\rho),$ with $a\in\mathbb{D},$ and we choose $0<r<1$ satisfying $\rho=d_h(0,r).$ Then the analytic function
		$\psi(z)=h(\varphi_{z_0}(rz))$ is univalent in $\mathbb{D}$ and
		\[S_{\psi}(z)=S_h(\varphi_{z_0}(rz))(\varphi_{z_0}'(rz))^2r^2,\quad z\in\mathbb{D},\]
		whence
		\[|S_{\psi}(0)|=|S_h(z_0)|(1-|z_0|^2)^2r^2>\delta r^2.\]
		We conclude from Lemma\,4 in \cite{Pau Pelaez} that there is $0<\widetilde{c}=\widetilde{c}(\delta,r)<1$ such that
		\[(1-|z|^2)^2|S_{\psi}(z)|>\frac{\delta r^2}{32},\]
		for all $z\in B(0, \widetilde{c}),$ and therefore
		\[(1-|z|^2)^2|S_h(\varphi_{z_0}(rz))||\varphi_{z_0}'(rz)|^2r^2>\frac{\delta r^2}{32},\]
		for all $z\in B(0, \widetilde{c}).$ Thus, by properties of the automorphisms of the disk we have
		\[(1-|\varphi_{z_0}(rz)|^2)^2||S_h(\varphi_{z_0}(rz))|>\frac{\delta}{32},\]
		for all $z\in B(0, \widetilde{c}).$ In consequence, by defining $\zeta=\varphi_{z_0}(rz),$ one has that
		\[(1-|\zeta|^2)^2||S_h(\zeta)|>\frac{\delta}{32},\]
		for all $\zeta\in D(z_0,\widetilde{\rho}),$ where $\widetilde{\rho}=d_h(0,r\widetilde{c}).$ The lemma follows having into account that the hyperbolic disk $D(z_0,\widetilde{\rho})$ contains the euclidean ball $B(z_0,r\widetilde{c}(1-|z_0|^2)/4).$ 
	\end{proof}
	\begin{theorem}\label{thm Qp spaces}
		Let $0<p<\infty$ and $f:\mathbb{D}\to\mathbb{C}$ be a uniformly locally univalent harmonic map with dilatation $w\in Q_p$ satisfying $\|w\|<1.$ Then $F=\log J_f\in \widetilde{\mathcal{Q}}_p$ if and only if $\mu=|S_f(z)|^2(1-|z|^2)^{2+p}$ is a $p-$Carleson measure.
	\end{theorem}
	\begin{proof}
		We suppose that \eqref{eq. def of Qs} holds, or equivalently,
		\begin{equation}\label{eq. aux thm Qs and carleson}
			\sup_{a\in\mathbb{D}}\int_{\mathbb{D}}\left|\frac{h''(z)}{h'(z)}-\frac{w'(z)\overline{w(z)}}{1-|w(z)|^2} \right|^2(1-|\varphi_a(z)|^2)^{p}dA(z)<\infty.
		\end{equation}
		By \eqref{eq. p carleson measure} it is sufficient prove that
		\begin{equation}\label{eq. conclusion}
			\sup_{a\in\mathbb{D}}\int_{\mathbb{D}}|S_f(z)|^2(1-|z|^2)^{2+p}|\varphi_a'(z)|^pdA(z)<\infty.
		\end{equation}
		If we denote by $A(f,a)$ the integral
		\[\int_{\mathbb{D}}|S_f(z)|^2(1-|z|^2)^{2+p}|\varphi_a'(z)|^pdA(z),\]
		we obtain, from \eqref{Schwarzian of f} and the Schwarz-Pick inequality, that
		\begin{equation}\label{eq. bound A(f,a)}
			A(f,a)\lesssim \int_{\mathbb{D}}\left(|F_{zz}(z)|^2+|F_z(z)|^4 \right)(1-|z|^2)^2(1-|\varphi_a(z)|^2)^pdA(z). 
		\end{equation}
		Now, on the one hand,
		\begin{multline*}
			\int_{\mathbb{D}}|F_z(z)|^4(1-|z|^2)^2(1-|\varphi_a(z)|^2)^pdA(z)=\\\int_{\mathbb{D}}[(1-|z|^2)|F_z(z)|]^2|F_z(z)|^2(1-|\varphi_a(z)|^2)^pdA(z),
		\end{multline*}
		whence
		\begin{equation*}
			\int_{\mathbb{D}}|F_z(z)|^4(1-|z|^2)^2(1-|\varphi_a(z)|^2)^pdA(z)\leq \|P_f\|^2\int_{\mathbb{D}}|F_z(z)|^2(1-|\varphi_a(z)|^2)^pdA(z).
		\end{equation*}
		We conclude from \eqref{eq. aux thm Qs and carleson} and the definition of $F_z$ that
		\begin{equation}\label{eq. bound int Fz}
			\sup_{a\in\mathbb{D}}\int_{\mathbb{D}}|F_z(z)|^4(1-|z|^2)^2(1-|\varphi_a(z)|^2)^pdA(z)<\infty. 
		\end{equation}
		On the other hand we show that
		\begin{equation}\label{eq. bound int Fzz}
			\sup_{a\in\mathbb{D}}\int_{\mathbb{D}}|F_{zz}(z)|^2(1-|z|^2)^2(1-|\varphi_a(z)|^2)^pdA(z)<\infty,
		\end{equation}
		which is equivalent to prove
		\begin{equation*}
			\sup_{a\in\mathbb{D}}\int_{\mathbb{D}}\left|  \left( \frac{h''(z)}{h'(z)}\right)'-\frac{w''(z)\overline{w(z)}}{1-|w(z)|^2}-\left(\frac{w'(z)\overline{w(z)}}{1-|w(z)|^2} \right)^2\right|^2(1-|z|^2)^2(1-|\varphi_a(z)|^2)^pdA(z)<\infty.
		\end{equation*}
		To this end we first note that 
		\[F=\log|h'|^2+\log(1-|w|^2)\qquad\text{and}\qquad F,w\in \widetilde{\mathcal{Q}}_p \]
		imply that $\log h'\in Q_p.$ Therefore, \eqref{eq. f in Qp iff ...} yields 
		\begin{equation}\label{eq. bound 1 for Sf}
			\sup_{a\in\mathbb{D}}\int_{\mathbb{D}}\left|\left( \frac{h''(z)}{h'(z)}\right) '\right|^2 (1-|z|^2)^2(1-|\varphi_a(z)|^2)^pdA(z)<\infty.
		\end{equation}
		Next, taking into account the conditions $w\in Q_p,$ $\|w\|<1$ and applying again \eqref{eq. f in Qp iff ...}, we see that
		\begin{equation}\label{eq. bound 2 for Sf}
			\sup_{a\in\mathbb{D}}\int_{\mathbb{D}}\left|\frac{w''(z)\overline{w(z)}}{1-|w(z)|^2}\right|^2 (1-|z|^2)^2(1-|\varphi_a(z)|^2)^pdA(z)<\infty.
		\end{equation}
		Now we use the Schwarz-Pick inequality to conclude that
		\begin{equation*}
			\int_{\mathbb{D}}\left|\frac{w'(z)\overline{w(z)}}{1-|w(z)|^2}\right|^4 (1-|z|^2)^2(1-|\varphi_a(z)|^2)^pdA(z)\leq\int_{\mathbb{D}}|w'(z)|^2 (1-|\varphi_a(z)|^2)^pdA(z),
		\end{equation*}
		which gives
		\begin{equation}\label{eq. bound 3 for Sf}
			\sup_{a\in\mathbb{D}}\int_{\mathbb{D}}\left|\frac{w'(z)\overline{w(z)}}{1-|w(z)|^2}\right|^4 (1-|z|^2)^2(1-|\varphi_a(z)|^2)^pdA(z)<\infty,
		\end{equation}
		since $w\in Q_p.$ From \eqref{eq. bound 1 for Sf}, \eqref{eq. bound 2 for Sf}, and \eqref{eq. bound 3 for Sf} it follows \eqref{eq. bound int Fzz}. Hence and \eqref{eq. bound int Fz} we obtain that
		\[\sup_{a\in\mathbb{D}}A(f,a)<\infty,\]
		and consequently $\mu=|S_f(z)|^2(1-|z|^2)^{2+p}$ is a $p-$Carleson measure.
		
		In order to prove the converse we will first use
		\begin{equation*}
			S_h=S_f-\frac{\overline{w}}{1-|w|^2}\left(w'\frac{h''}{h'}-w'' \right)
			+\frac{3}{2}\left(\frac{w'\overline{w}}{1-|w|^2} \right)^2
		\end{equation*}
		to prove that $|S_h(z)|^2(1-|z|^2)^{2+p}$ is a $p-$Carleson measure.
		
		In view of $\|P_h\|<\infty,$ which is a consequence of the fact that $f$ is uniformly locally univalent, and the condition $\|w\|<1,$ we can obtain
		\begin{multline*}
			\int_{\mathbb{D}}\left| \frac{\overline{w(z)}}{1-|w(z)|^2}\right|^2\left|w'(z)\frac{h''(z)}{h'(z)} \right|^2(1-|z|^2)^{2+p}|\varphi_a'(z)|^pdA(z)\\
			\lesssim \int_{\mathbb{D}}|w'(z)|^2(1-|z|^2)^{p}|\varphi_a'(z)|^pdA(z),
		\end{multline*}
		whence in virtue of $w\in Q_p,$ we conclude that
		\begin{equation}\label{eq. bound 1 for Sh}
			\sup_{a\in\mathbb{D}}\int_{\mathbb{D}}\left| \frac{\overline{w(z)}}{1-|w(z)|^2}\right|^2\left|w'(z)\frac{h''(z)}{h'(z)} \right|^2(1-|z|^2)^{2+p}|\varphi_a'(z)|^pdA(z)<\infty.
		\end{equation}
		Now we use again $\|w\|<1$ to obtain
		\[\int_{\mathbb{D}}\left| \frac{\overline{w(z)}}{1-|w(z)|^2}\right|^2\left|w''(z) \right|^2(1-|z|^2)^{2+p}|\varphi_a'(z)|^pdA(z)\lesssim \int_{\mathbb{D}}\left|w''(z) \right|^2(1-|z|^2)^{2+p}|\varphi_a'(z)|^pdA(z),\]
		which, because of \eqref{eq. f in Qp iff ...} and $w\in  Q_p,$ gives
		\begin{equation}\label{eq. bound 2 for Sh}
			\sup_{a\in\mathbb{D}}\int_{\mathbb{D}}\left| \frac{\overline{w(z)}}{1-|w(z)|^2}\right|^2\left|w''(z) \right|^2(1-|z|^2)^{2+p}|\varphi_a'(z)|^pdA(z)<\infty.
		\end{equation}
		It follows from \eqref{eq. bound 3 for Sf},\eqref{eq. bound 1 for Sh}, and \eqref{eq. bound 2 for Sh} that $|S_h(z)|^2(1-|z|^2)^{2+p}$ is a $p-$Carleson measure. Thus Proposition\,\ref{Prop Thm of Pau-Pelaez for locally univalent} implies that $\log h'\in Q_p,$ whence in virtue of the equality $F=\log|h'|^2+\log(1-|w|^2)$ and the hypothesis $w\in Q_p,$ it follows that $F=\log J_f\in \widetilde{\mathcal{Q}}_{p}.$ 
	\end{proof}
	With some minor changes in the proof of the first implication of Theorem\,\ref{thm Qp spaces}, we can extend to the setting of harmonic mappings the implication ``only if" of part ii) of Theorem\,A. In this case, instead of \eqref{eq. f in Qp iff ...} we use the fact that for $n\geq 1$ and $0<p<\infty,$ $h\in Q_{p,0}$ if and only if
	\begin{equation}\label{eq. f in Qp0 iff ...}
		\lim_{|a|\to 1}\int_{\mathbb{D}}|h^{n}(z)|^2\left(1-|\varphi_a(z)|^2 \right)^p(1-|z|^2)^{2n-2}dA(z)= 0,
	\end{equation}
	see \cite[Theorem 2]{Qs and n derivative}. The final result would be the following:
	\begin{theorem}\label{thm Qp0 spaces}
		Let $0<p<\infty$ and $f:\mathbb{D}\to\mathbb{C}$ be a uniformly locally univalent harmonic map with dilatation $w\in Q_{p,0}$ satisfying $\|w\|<1.$ If $F=\log J_f\in \widetilde{\mathcal{Q}}_{p,0},$ then $\mu=|S_f(z)|^2(1-|z|^2)^{2+p}$ is a vanishing $p-$Carleson measure.
	\end{theorem}
	\begin{remark}
		One can prove that the harmonic function $f$ constructed in Example\,\ref{exa shear construction} satisfies the hypotheses of the above theorem for $0<p<1,$ therefore $\mu=|S_f(z)|^2(1-|z|^2)^{2+p}$ is a vanishing $p-$Carleson measure, if $0<p<1.$
	\end{remark}
	A harmonic version of the other implication of part ii) of Theorem\,A can be obtained by using the method employed in Theorem\,\ref{Thm 2}.
	\begin{theorem}\label{Thm 5}
		Let $0<p\leq 1$ and $f=h+\overline{g}:\mathbb{D}\to\mathbb{C}$ be a univalent harmonic map such that its dilatation $w\in \mathcal{Q}_{p,0}$ satisfies $\|w\|<1$ and assume that $\varphi:=\varphi_{\lambda}=h+\lambda g$ is a conformal map with $\varphi(\partial\mathbb{D})$ a closed Jordan curve for some $\lambda\in\partial\mathbb{D}.$ If $|S_f(z)|^2(1-|z|^2)^{2+p}dA(z)$ is a vanishing $p-$Carleson measure, then $F=\log J_f\in \widetilde{\mathcal{Q}}_{p,0}.$ 
	\end{theorem}
	\begin{proof}
		We can verify by a straightforward calculation, using \eqref{Schwarzian of f}, that
		\[d\mu_f:=|S_f(z)|^2(1-|z|^2)^{2+p}dA(z)\]
		is a vanishing $p-$Carleson measure if and only if
		\[d\mu_h:=|S_h(z)|^2(1-|z|^2)^{2+p}dA(z)\]
		is a vanishing $p-$Carleson measure, which is equivalent to the fact that
		\[d\mu_{\varphi}:=|S_\varphi(z)|^2(1-|z|^2)^{2+p}dA(z)\]
		is a vanishing $p-$Carleson measure, by  \eqref{eq. Sh and S(phi)}. In both cases \eqref{eq. f in Qp0 iff ...} is applied. Then part ii) of Theorem\,A implies that $\log\varphi'\in\mathcal{Q}_{p,0},$ whence
		\[\log h'= \log\varphi'-\log(1+\lambda w)\in \mathcal{Q}_{p,0}. \]
		The conclusion of the theorem follows from
		\[|F_z|^2\lesssim \left| \frac{h''}{h'} \right|^2 + \left( \frac{|w'\overline{w}|}{1-|w|^2} \right)^2\]
		and the conditions $w\in \mathcal{Q}_{p,0}$ and $\|w\|<1.$
	\end{proof}
	
	\section{On harmonic Besov-type and $\widetilde{\mathcal{Q}}_p-$type mappings class}\label{section results in harmonic spaces}
	Our main objective in this section is to introduce the Besov-type and $\widetilde{\mathcal{Q}}_p-$type classes of harmonic mappings and study their connection with the Bloch-type class of harmonic mappings defined in \cite{Bloch type map}, see also \cite{LP v bloch and v bloch type} for a generalization. A smooth function $F:\mathbb{D}\to\mathbb{C}$ is said to be Bloch-type if
	\begin{equation}\label{eq. def H type bloch }
		\beta_2(F):=\sup_{z\in\mathbb{D}}(1-|z|^2)\sqrt{|J_F(z)|\;}<\infty.
	\end{equation}
	We use $\mathcal{BT}$ to denote this class of functions, which clearly contains the analytic Bloch space $\mathcal{B}$ and, although it is not a linear space, it also contains the Bloch space $\widetilde{\mathcal{B}},$ being this fact an immediate consequence of $J_F^{1/2}\leq |F_{z}|+|F_{\overline{z}}|.$ Example\,2 in \cite{Bloch type map} shows that there are harmonic functions in $\mathcal{BT}$ that do not belong to the space $\widetilde{\mathcal{B}}.$ It is natural to define $\mathcal{BT}_0$ as the set of all smooth function $F:\mathbb{D}\to\mathbb{C}$ such that \[\lim_{|z|\to 1^{-}}(1-|z|^2)|J_F(z)|^{1/2}=0.\]
	Following these ideas, we define the class of $\widetilde{\mathcal{Q}}_p-$type functions, $0<p<\infty,$ denoted by $\mathcal{QT}_p,$ as the class of smooth functions $F:\mathbb{D}\to\mathbb{C}$ such that
	\begin{equation*}\label{eq. def Qs type}
		\sup_{a\in\mathbb{D}}\int_{\mathbb{D}}|J_F(z)|(1-|\varphi_a(z)|^2)^pdA(z)<\infty.
	\end{equation*}
	It is clear that $Q_p\subseteq \widetilde{\mathcal{Q}}_p\subseteq \mathcal{QT}_p,$ for all $0<p<\infty.$ Similarly, we will say that $F$ belongs to the class of Besov-type functions $\mathcal{BT}_p,$ $1<p<\infty,$ if
	\begin{equation}\label{eq. def Bp type}
		\beta(F):=\int_{\mathbb{D}}|J_F(z)|^{p/2}(1-|z|^2)^{p-2}dA(z)<\infty.
	\end{equation}
	Here also it is easy to see that $\mathcal{B}_p\subseteq \widetilde{\mathcal{B}}_p\subseteq \mathcal{BT}_p.$
	
	The following proposition exhibits some important properties of the previous classes.
	\begin{proposition}\label{prop affine and linear invariance}
		The classes $\mathcal{BT}_p$, $p>1,$ and $\mathcal{QT}_p,$ $p>0,$ are affine and linearly invariant.
	\end{proposition}
	\begin{proof} For all automorphism $\sigma$ of $\mathbb{D},$ $|J_{f\circ \sigma}(z)|=|J_f(\sigma(z))||\sigma'(z)|^2,$ $z\in\mathbb{D}.$ Hence,
		\begin{align*}
			\int_{\mathbb{D}}|J_{F\circ\sigma}(z)|^{p/2}(1-|z|^2)^{p-2}dA(z)&=\int_{\mathbb{D}}|J_F(\sigma(z))|^{p/2}|\sigma'(z)|^p(1-|z|^2)^{p-2}dA(z)\\
			&=\int_{\mathbb{D}}|J_F(\sigma(z))|^{p/2}(1-|\sigma(z)|^2)^{p-2}|\sigma'(z)|^2dA(z)\\
			&=\int_{\mathbb{D}}|J_F(z)|^{p/2}(1-|z|^2)^{p-2}dA(z),
		\end{align*}
		which shows that the class $\mathcal{BT}_p$ is linearly invariant. This  is, $F\circ\sigma\in \mathcal{BT}_p ,$ for all automorphism $\sigma$ of $\mathbb{D}$ and $F\in\mathcal{BT}_p.$ Also it is easy to see the affine invariance of $\mathcal{BT}_p,$ which means that $aF+b\overline{F}\in\mathcal{BT}_p, $ for all $a,b\in\mathbb{C}$ and $F\in\mathcal{BT}_p. $ This follows from the equality $J_{aF+b\overline{F}}=(|a|^2-|b|^2)J_F.$ To show the second statement we observe that
		\begin{align*}
			\int_{\mathbb{D}}|J_{F\circ\sigma}(z)|(1-|\varphi_a(z)|^2)^pdA(z)&=\int_{\mathbb{D}}|J_F(\sigma(z))|(1-|\varphi_a(z)|^2)^p|\sigma'(z)|^2dA(z)\\
			&=\int_{\mathbb{D}}|J_F(\zeta)|(1-|\varphi_a(\sigma^{-1}(\zeta))|^2)^pdA(\zeta)\\
			&=\int_{\mathbb{D}}|J_F(z)|(1-|\varphi_b(z)|^2)^pdA(z),
		\end{align*}
		for some $b:=b(a,\sigma)\in\mathbb{D},$ which implies that $\mathcal{QT}_p$ is linearly invariant. The affine invariance of $\mathcal{QT}_p$ follows also of the equality $J_{aF+b\overline{F}}=(|a|^2-|b|^2)J_F.$ 
	\end{proof}
	We remark that the affine and linear invariance of $\mathcal{BT}$ was proven in \cite{Bloch type map}. Those properties allow us to prove the following proposition, which generalizes a statement in \cite{Book_Pommerenke} page 73.
	
	\begin{proposition}  
		Let $f$ be a harmonic mapping in $\in\mathcal{BT}$ and $\|w\|<1,$ then 
		\[|f(z_1)-f(z_2)|\leq \left(\frac{1+\|w\|}{1-\|w\|}\right)^{1/2}\beta_2(f)d_h(z_1,z_2),\]
		for all $z_1,z_2\in \mathbb{D}.$
	\end{proposition}
	\begin{proof}
		Since $\beta_2(f)=\beta_2(f\circ \sigma)$ for $\sigma \in Aut(\mathbb D)$, we can assume that $z_1=z$ and $z_2=0.$ Thus,  
		\begin{align*}
			|f(z)-f(0)|&\leq \int_{0}^{1}\left|(f_z(tz)z+f_{\overline{z}}(tz)\overline{z})\right|dt\\
			&\leq\int_{0}^{1}|f_z(tz)|(1-|w(tz)|^2)^{1/2}\frac{1+|w(tz)|}{(1-|w(tz)|^2)^{1/2}}|z|dt\\
			&=\int_{0}^{1}|J_f^{1/2}(tz)|\left(\frac{1+|w(tz)|}{1-|w(tz)|}\right)^{1/2}|z|dt\\
			&\leq \left(\frac{1+\|w\|}{1-\|w\|}\right)^{1/2}\beta_2(f)\int_{0}^{1} \frac{|z|}{1-t^2|z|^2}dt.
		\end{align*}
		After integrating we obtain
		\begin{align*}
			|f(z)-f(0)|&\leq \frac{1}{2}\left(\frac{1+\|w\|}{1-\|w\|}\right)^{1/2}\beta_2(f)\ln\left(\frac{1+|z|}{1-|z|}\right)\\
			&=\left(\frac{1+\|w\|}{1-\|w\|}\right)^{1/2}\beta_2(f)d_h(0,z),
		\end{align*}
		which ends the proof.
	\end{proof}
	
	\begin{example}
		Let $f=h+\overline{g}$ with $h(z)=2(1-z)^{-1/2}$ and $w(z)=z$, $z\in \mathbb{D}.$ We observe that 
		\[J_f(z)=|h'(z)|^2(1-|z|^2)=\frac{1-|z|^2}{|1-z|^3}.\]
		\begin{itemize}
			\item [i)] $f\in \mathcal{BT}$ since
			\[(1-|z|^2)J_f(z)^{1/2}=\left(\frac{1-|z|^2}{|1-z|}\right)^{3/2}\leq2^{3/2},\]
			for all $z\in\mathbb{D}.$
			\item [ii)] $f\not\in \mathcal{BT}_0.$ In effect, if $z=x\in(-1,1),$
			\[\lim_{x\to 1^-}(1-x^2)J_f(x)^{1/2}=2^{3/2}.\]
			\item [iii)] For $p>1,$ $f\in \mathcal{QT}_p.$ This follows from the fact that
			\begin{align*}
				\int_{\mathbb{D}}J_f(z)(1-|\varphi_a(z)|^2)^pdA&=\int_{\mathbb{D}}\frac{1-|z|^2}{|1-z|^3}(1-|\varphi_a(z)|^2)^pdA\\
				&=\int_{\mathbb{D}}\left(\frac{1-|z|^2}{|1-z|}\right)^{3}(1-|\varphi_a(z)|^2)^{p-2}|\varphi'_a(z)|^2dA\\
				&\leq8\int_{\mathbb{D}}(1-|\zeta|^2)^{p-2}dA(\zeta)<\infty, 
			\end{align*}
			for all $a\in\mathbb{D},$ where $\zeta=\varphi_a(z).$
			\item [iv)] In part ii) of Remark\,\ref{rem 4} below we will prove that $f\not\in \mathcal{BT}_p,$ $p>1.$ 
		\end{itemize}
	\end{example}
	\begin{remark}\label{rem 4}
		The previous example suggests that we should study how the classes defined above are related.
		\begin{itemize}
			\item [i)] For $p>1,$ $\mathcal{BT}\subseteq \mathcal{QT}_p.$ In effect, if $F\in \mathcal{BT},$ for all $a\in\mathbb{D}$ \eqref{eq. def H type bloch } implies
			\begin{align*}
				\int_{\mathbb{D}}|J_F(z)|(1-|\varphi_a(z)|^2)^pdA(z)&=\int_{\mathbb{D}}|J_F(z)^{1/2}(1-|z|^2)|^2(1-|z|^2)^{p-2}|\varphi'_a(z)|^pdA(z)\\
				&\leq \beta_2(F)^2\int_{\mathbb{D}}(1-|z|^2)^{p-2}|\varphi'_a(z)|^pdA(z)\\
				&=\beta_2(F)^2\int_{\mathbb{D}}(1-|\varphi_a(z)|^2)^{p-2}|\varphi'_a(z)|^2dA(z).
			\end{align*}
			Hence, the change of variables $\zeta=\varphi_a(z)$ gives
			\[\int_{\mathbb{D}}|J_F(z)|(1-|\varphi_a(z)|^2)^pdA(z)\leq \beta_2(F)^2\int_{\mathbb{D}}(1-|\zeta|^2)^{p-2}dA(\zeta)=C(F),\]
			where $C(F)$ denotes a constant depending only on $F.$ Thus, $F \in \mathcal{QT}_p$, for all $p>1.$
			
			\item [ii)] If $f=h+\overline{g}\in\mathcal{BT}_p$, $p>1$, and $w(z)=z$, then $f\in \mathcal{BT}_0.$ Indeed, $f\in \mathcal{BT}_p$ implies
			\[\int_{\mathbb{D}}[(1-|z|^2)^{3/2}|h'(z)|]^p(1-|z|^2)^{-2}dA(z)<\infty.\] So, if $k':=(h')^{2/3}$, then $k\in \mathcal{B}_{3p/2}$ and therefore $k\in \mathcal{B}_0$.      
			In consequence,
			\begin{align*}
				\lim_{|z|\to 1^{-}}(1-|z|^2)J_f(z)^{1/2}&=\lim_{|z|\to 1^{-}}[(1-|z|^2)|h'(z)|^{2/3}]^{3/2}\\
				&=\lim_{|z|\to 1^{-}}[(1-|z|^2)|k'(z)|]^{3/2}=0.
			\end{align*}
			Hence, $f\in \mathcal{BT}_0.$
			
			The above statement remains valid if $(1-|w|^2)\,\thicksim\,(1-|z|^2).$ In particular, 
			\begin{itemize}
				\item [a)] If $w(z)=z^n, n\in \mathbb {N}$, it is sufficient to observe that 
				\[1- |z^n|^2=(1-|z|^2)(1+|z|^2+|z|^4+\cdots+|z|^{2n-2})\geq 1-|z|^2,\]
				and we proceed as above. 
				\item [b)] If $w(z)=(z+a)/(1+\overline{a}z),$ we consider $F=f-\overline{af}$. By Proposition\,\ref{prop affine and linear invariance}, if $f\in \mathcal{BT}_p$ then $F\in \mathcal{BT}_p$. Since $w_F=Id,$ Part ii) above implies that $F\in \mathcal{BT}_0$ and so $f\in \mathcal{BT}_0.$
			\end{itemize}
		\end{itemize}
	\end{remark}
	
	\begin{example}
		The authors in \cite{ACS} introduce the class $NH_{\mu}^{0},$ $ 0<\mu\leq1,$ of locally univalent sense preserving harmonic mappings $f$, defined in the unit disk,  for which 
		\[|S_f(z)|+\frac{|w'(z)|^2}{(1-|w(z)|^2)^2}\leq \frac{2\mu}{(1-|z|^2)^2},\]
		with $\nabla J_f(0,0)=(0,0)$. By Using Theorem\,2.2 in \cite{ACS}, for $f\in NH_{\mu}^{0}$ with $J_f(0)=1$ we have that
		\begin{itemize}
			\item [i)] If $\mu=1$ then 
			\[1-|z|^2\leq J_f^{1/2}(z)\leq(1-|z|^2)^{-1},\]
			whence $f\in\mathcal{BT}$ and even more $\beta_2(f)=1.$
			\item [ii)] If $0<\mu<1$ then 
			\begin{equation}\label{ex}
				\frac{[(1+|z|)^{\beta}+(1-|z|)^\beta]^2}{4(1-|z|^2)^{\beta-1}}\leq J_f^{1/2}(z)\leq\frac{4(1-|z|^2)^{\beta-1}}{[(1+|z|)^{\beta}+(1-|z|)^\beta]^2},
			\end{equation}
			where $\beta=\sqrt{1-\mu\,}.$ Thus,  $1\leq \beta_2(f)\leq 4.$
			\item [$(iii)$] If $0<\mu<1$ it follows from (\ref{ex}) that
			\begin{align*}
				\int_{\mathbb{D}}J_f^{p/2}(z)(1-|z|^2)^{p-2}dA&\leq 4^p\int_{\mathbb{D}}(1-|z|^2)^{p(\beta-1)}(1-|z|^2)^{p-2}dA\\
				&=4^p\int_{\mathbb{D}}(1-|z|^2)^{p\beta-2}dA,
			\end{align*}
			and it is known that the integral on the right side is finite if $p\beta-2>-1$. We conclude that $f\in\mathcal{BT}_p$ if $p > 1/\beta.$
		\end{itemize}
	\end{example}
	
	\textbf{Comments.} We believe that the topic introduced in Section\,\ref{section results in harmonic spaces} could serve as the foundation for a more extensive study in the future. For instance, the inclusion relationship $\mathcal{QT}_p \subseteq \mathcal{BT}$, where $p > 1$, though presumed, remains unconfirmed. Similarly, the containment of $\mathcal{BT}_p$, $p > 1$, in $\mathcal{BT}_0$, for any dilation $w$, is yet to be established.

\end{document}